\documentclass[11pt]{article}

\usepackage{cite}

\usepackage[english]{babel}
\usepackage{amsmath}
\usepackage{amsthm}
\usepackage{amssymb}
\usepackage{tikz}
\usepackage{tikz-cd}
\usepackage{csquotes}

\usepackage{hyperref}

\usepackage{geometry}
\usepackage{enumitem}
\newtheorem{thm}{Theorem}[section]
\newtheorem*{thm*}{Theorem}

\newtheorem{theorem}{Theorem}

\newtheorem{proposition}[thm]{Proposition}
\newtheorem{lem}[thm]{Lemma}
\newtheorem{cor}[thm]{Corollary}

\theoremstyle{definition}

\newtheorem{defn}[thm]{Definition}

\theoremstyle{remark}
\newtheorem{rem}[thm]{\textsc{Remark}}

\title{
Identifying the homotopy fiber of a map of semi-Segal spaces
}

\author{Yuxun Sun}
\date{October 24, 2024}

\begin{document}

\maketitle
 
\begin{abstract}
We prove a version of Quillen's theorems for a map of semi-Segal spaces. We construct a bi-semi-simplicial resolution similar to the one associated to a functor of non-unital topological categories. As a consequence we can represent the homotopy fiber of a map of semi-Segal spaces as the geometric realization of a certain semi-simplicial space. 
\end{abstract}


\section{Introduction}
Let $F : C \rightarrow D$ be a functor between two small categories. We can associate to this functor a continuous map $\textbf{B}F: \textbf{B} C \rightarrow \textbf{B} D$ between the classifying spaces of $C$ and $D$ in the sense of \cite[pp. 106]{segal1968classifying}. One natural question to ask is when $F$ induces a homotopy equivalence between the classifying spaces. Another interesting question is under what conditions we can identify the homotopy fiber of the map $\textbf{B}F$ with the classifying space of some category. Quillen answered both of these questions by proving what are later referred to as Quillen's Theorem A and Theorem B \cite[pp. 93-98]{quillen2006higher}. In particular, he showed that under some conditions, the homotopy fiber over some object $d \in D$ can be identified with $\textbf{B} F/d$, where $F/d$ is the comma category over $d$. 

Since then, Quillen's theorems have been successfully generalized in various settings. In \cite{roberts2024homotopy}, the author extended Quillen’s Theorem A to apply to a continuous functor between (unital) topological categories. In \cite[Theorem 4.7 and 4.9]{ebert2019semisimplicial}, the authors proved another version of Quillen's Theorem A and Theorem B for a functor of non-unital topological categories. In this case, the classifying space $\textbf{B} \mathcal{C}$ of a non-unital topological category $\mathcal{C}$ is the fat geometric realization of the semi-simplicial nerve $N_\bullet \mathcal{C}$ \cite[Definition 3.3]{ebert2019semisimplicial}. 

The goal of this paper is to prove a version of Quillen's theorems for a map of semi-Segal spaces. The notion of semi-Segal spaces is a natural generalization of the semi-simplicial nerve of non-unital topological categories. Recall that a \emph{semi-simplicial} space is a functor $\Delta_{\text{inj}}^{\text{op}} \rightarrow \textbf{Top}$, where $\Delta_{\text{inj}} \subset \Delta$ is the subcategory of the simplex category consisting of only injective maps. A \emph{semi-Segal} space $C_\bullet$ is a semi-simplicial space which satisfies the Segal condition: for each $n \geq 2$, the canonical map
    \begin{center}
        $C_n \rightarrow C_1 \times_{C_0}^h ... \times_{C_0}^h C_1$
    \end{center}
is a weak homotopy equivalence. Here $C_1 \times_{C_0}^h ... \times_{C_0}^h C_1$ is the 
n-fold iterated homotopy fiber product. For example, the semi-simplicial nerve $N_\bullet \mathcal{C}$ of a non-unital \emph{locally fibrant} (\cite[Definition 2.1]{steimle2021additivity}) topological category $\mathcal{C}$ is a semi-Segal space. 

By a \emph{map of semi-Segal spaces}, we will mean a map of the underlying semi-simplicial spaces. The study of maps of semi-Segal spaces has important applications in cobordism categories, which are the topological categories originally studied in \cite{galatius2009homotopy}. These cobordism categories can be viewed and studied as semi-Segal spaces. In \cite[Theorem 2.11]{steimle2021additivity}, Steimle formulated a local-to-global principle for semi-Segal spaces, which has proven to be a powerful tool for obtaining homotopy cartesian diagrams of classifying spaces of cobordism categories. Our results in this paper complement those of Steimle by providing a systematic approach to identifying the homotopy fiber of a map between semi-Segal spaces. On the other hand, our results generalize \cite[Theorem 4.7 and 4.9]{ebert2019semisimplicial}. As we will see in Section \ref{sec 3}, our construction of the bi-semi-simplicial resolution differs from that in \cite[Section 4]{ebert2019semisimplicial}, as there is no strict composition of two 1-simplices in a semi-Segal space.

Given a semi-Segal space $C_\bullet$, we denote by $\Vert C_\bullet \Vert$ the geometric realization of $C_\bullet$ as a semi-simplicial space in the sense of \cite[pp. 4]{ebert2019semisimplicial}. We will prove the following two theorems:

\begin{theorem}   \label{thm A}
      Let $F : C_\bullet \rightarrow D_\bullet$ be a map of semi-Segal spaces. Suppose that 
    \begin{itemize}
        \item $D_\bullet$ is weakly unital; 
        \item $\Vert (F/^hy_0)_{\bullet} \Vert$ is contractible for each $y_0 \in D_0$.
    \end{itemize}
    Then the map $\Vert F \Vert : \Vert C_\bullet \Vert \rightarrow \Vert D_\bullet \Vert$ is a weak homotopy equivalence.
\end{theorem}

\begin{theorem}   \label{thm B}
     Let $F : C_\bullet \rightarrow D_\bullet$ be a map of semi-Segal spaces. Suppose that 
    \begin{itemize}
        \item $D_\bullet$ is weakly unital; 
        \item for every 1-simplex $y \in D_1$, the induced map $\Vert (F/^hy)_\bullet \Vert  \rightarrow  \Vert (F/^hd_0y)_\bullet \Vert$ is a weak homotopy equivalence.
    \end{itemize}
    Then for every $y_0 \in D_0$, there exists a homotopy fibration sequence 
    \begin{center}
        $ \Vert (F/^hy_0)_\bullet \Vert \longrightarrow \Vert C_\bullet \Vert \overset{\Vert F \Vert}{\longrightarrow} \Vert D_\bullet \Vert$
    \end{center} 
\end{theorem}

We also show that the assumptions in Theorem \ref{thm B} can be satisfied in a special case, as expressed in the next corollary. 

\begin{cor}    \label{cor c}
    Let $C_\bullet \rightarrow D_\bullet$ be a map of semi-Segal spaces. Suppose that $D_\bullet$ is a weakly unital semi-Segal space in which every 1-simplex in $D_1$ is cartesian. Then for every $y_0 \in D_0$, there exists a homotopy fibration sequence 
    \begin{center}
        $ \Vert (F/^hy_0)_\bullet \Vert \longrightarrow \Vert C_\bullet \Vert \overset{\Vert F \Vert}{\longrightarrow} \Vert D_\bullet \Vert$
    \end{center} 
\end{cor}

\subsection*{Structure of the paper}
Our methods closely follow \cite{ebert2019semisimplicial} and \cite{steimle2021additivity}. In section \ref{sec 2}, we recall some necessary definitions and results for semi-Segal spaces. We briefly discuss Reedy fibrations and Reedy fibrant replacement. In section \ref{sec 3}, we construct a bi-semi-simplicial resolution similar to the one in \cite[Definition 4.1]{ebert2019semisimplicial}. This will serve as the main tool for this paper. Sections \ref{sec 4} and \ref{sec 5} will be devoted to proving Theorem \ref{thm A}, Theorem \ref{thm B} and Corollary \ref{cor c}.

\subsection*{Acknowledgments}
I would like to express my gratitude to my advisors, Bernard Badzioch and Wojciech Dorabiała, for their guidance and unwavering support throughout the writing process. I am also grateful for their numerous helpful comments and encouragement.

\section{Semi-Segal spaces}   \label{sec 2}
 We first review the notion of equivalences and weak unitality of semi-Segal spaces. Recall that a commutative square of topological spaces
\[
\begin{tikzcd}[column sep=10ex]
    A \arrow{d}{f} \arrow{r}{k_1} & D \arrow{d}{g}\\
    B \arrow{r}{k_0} & C
\end{tikzcd}
\]
is \emph{locally homotopy cartesian at}  $b_0 \in B$ if the induced map on the vertical homotopy fibers
\begin{center}
    $\text{hofib}_{b_0}(f) \rightarrow \text{hofib}_{k_0(b_0)} (g)$
\end{center}
is a weak homotopy equivalence. A square $(ABCD)$ is homotopy cartesian if and only if it is locally homotopy cartesian at every $b_0 \in B$. 

\begin{defn}   \label{def T-cartesian}
    Let $C_\bullet$ be a semi-Segal space. We say that a 1-simplex $f \in C_1$ is \emph{cartesian} if the square
    \[
    \begin{tikzcd}[column sep=10ex]
        C_2 \arrow{r}{d_1} \arrow{d}{d_0} & C_1 \arrow{d}{d_0} \\
        C_1 \arrow{r}{d_0} & C_0
    \end{tikzcd}
    \]
    is locally homotopy cartesian at $f$. Similarly, $f$ is said to be \emph{cocartesian} if 
    \[
    \begin{tikzcd}[column sep=10ex]
        C_2 \arrow{r}{d_1} \arrow{d}{d_2} & C_1 \arrow{d}{d_1} \\
        C_1 \arrow{r}{d_1} & C_0
    \end{tikzcd}
    \]
    is locally homotopy cartesian at $f$. If $f$ is both cartesian and cocartesian, then $f$ is called an \emph{equivalence} in $C_\bullet$. Denote by $C_1^\simeq \subset C_1$ the subspace of equivalences. 
\end{defn}

 Recall that from \cite[Definition 2.9]{steimle2021additivity}, a semi-Segal space $C_\bullet$ is called \emph{weakly unital} if the map
    \begin{center}
        $d_0 : C_1^{\simeq} \rightarrow C_0$
    \end{center}
is 0-connected. A map of semi-Segal spaces is called \emph{weakly unital} if it maps an equivalence to an equivalence. 

In the remainder of this section, we shall only assume that $C_\bullet$ is a semi-simplicial space (not necessarily semi-Segal). For any semi-simplicial set $A$, we define $C(A)$ to be the space 
\begin{center}
    $C(A) := \text{Nat} (A, C_\bullet)$
\end{center}
where $\text{Nat}(A, C_\bullet)$ is the set of all natural transformations between $A$ and $C_\bullet$, with the subspace topology from $\underset{n}{\prod} \,\,\text{Map}(A_n, X_n)$. Here $A_n$ is equipped with the discrete topology. Note that a map $g : A \rightarrow B$ of semi-simplicial sets induces a map of spaces $g^*: C(B) \rightarrow C(A)$, defined by pre-composition with $g$.

\begin{defn}    \label{def reedy fib}
    Let $F : C_\bullet \rightarrow D_\bullet$ be a map of semi-simplicial spaces. We say that $F$ is a \emph{Reedy fibration} if for any inclusion $A \subset B$ of semi-simplicial sets, the map
    \begin{center}
        $C(B) \rightarrow D(B) \times_{D(A)} C(A) $
    \end{center}
    is a Serre fibration. $C_\bullet$ is called \emph{Reedy fibrant} if the map $T : C_\bullet \rightarrow * $ is a Reedy fibration, where $*$ is the one-point semi-simplicial space.
\end{defn}

\begin{rem}   \label{rem reedy fib}
    If $C_\bullet$ is Reedy fibrant, then for all $0 \leq i \leq n$, the structure map $d_i : C_n \rightarrow C_{n-1}$ is a fibration. 
\end{rem}

Recall that by \cite[Lemma 3.8]{steimle2021additivity}, a semi-simplicial map $F : C_\bullet \rightarrow D_\bullet$ admits a certain factorization: $F$ fits into the following triangle 
\[
    \begin{tikzcd}
        & \widetilde{C}_\bullet   \arrow{dr}{\widetilde{F}} & \\
        C_\bullet \arrow{ur}{\simeq}  \arrow{rr}{F} & & D_\bullet
    \end{tikzcd}
\]
where $C_\bullet \rightarrow \widetilde{C}_\bullet$ is a level equivalence, and $\widetilde{F}$ is a Reedy fibration. So in particular, for every semi-simplicial space $C_\bullet$,  there exists a level equivalence $C_\bullet \rightarrow C_\bullet^{\text{Rf}}$ where $C_\bullet^{\text{Rf}}$ is Reedy fibrant.

\begin{lem}  \label{lemma stable under reedy}
    Let $C_\bullet \rightarrow C_\bullet^{\text{Rf}}$ be a level equivalence of semi-Segal spaces,  where $C_\bullet^{\text{Rf}}$ is Reedy fibrant. Then the following statements are true:
    \begin{enumerate}[label=(\roman*)]
        \item If $C_\bullet$ is weakly unital, then so is $C_\bullet^{\text{Rf}}$; 
        \item If every 1-simplex in $C_1$ is cartesian(cocartesian), then every 1-simplex in $C_1^{\text{Rf}}$ is also cartesian(cocartesian). 
    \end{enumerate}
\end{lem}

\proof
Part (i) follows immediately from \cite[Lemma 3.9]{steimle2021additivity}. We prove part (ii) here. Consider the following cube
\[
\begin{tikzcd}
    & C_2 \arrow{rr}{d_1} \arrow{dl} \arrow{dd}[near end]{d_0} & & C_1 \arrow{dd}{d_0} \arrow{dl} \\
    C_2^{\text{Rf}} \arrow{rr}[near start]{d_1} \arrow{dd}{d_0} & & C_1^{\text{Rf}} \arrow{dd}[near start]{d_0} & \\
    & C_1 \arrow{rr}[near start]{d_0} \arrow{dl} & & C_0 \arrow{dl} \\
    C_1^{\text{Rf}} \arrow{rr}{d_0} & & C_0^{\text{Rf}} & 
\end{tikzcd}
\]
where the back face is homotopy cartesian. The left and the right faces are also homotopy cartesian, since $C_\bullet \rightarrow C_\bullet^{\text{Rf}}$ is a level equivalence. Therefore, by the 2-out-of-3 property (for example, see \cite[Lemma 3.3]{steinebrunner2022locally}), and the fact that $C_1 \rightarrow C_1^{\text{Rf}}$ is 0-connected, then the front face is also homotopy cartesian. The other case is proved in a similar way. 
\qed 

\section{A bi-semi-simplicial resolution}   \label{sec 3}

Let $F : C_\bullet \rightarrow D_\bullet$ be a map of semi-simplicial spaces. For each $p, q \geq 0$, we define $(F/^hD)_{p,q}$ to be the homotopy fiber product of the diagram
\begin{equation}   \label{dg 6.11}
    C_p \overset{F_p}{\longrightarrow} D_p \overset{d_{p+1}...d_{p+q+1}}{\xleftarrow{\hspace{1cm}}} D_{p+q+1}
\end{equation}
An element in $(F/^hD)_{p,q}$ is of the form $(x, \gamma, y) $, where $x \in C_p$, $y \in D_{p+q+1}$, and  $\gamma$ is a path in $D_p$ with $\gamma (0) = F(x)$ and $\gamma (1) = d_{p+1} ... d_{p+q+1}(y)$. 

We show that $(F/^hD)_{\bullet,\bullet}$ forms a bi-semi-simplicial space. In the $p$-direction, for all $0 \leq i \leq p$, the map $d_i^{(p)} : (F/^hD)_{p,q} \rightarrow (F/^hD)_{p-1,q}$ is induced from the diagram
\[
\begin{tikzcd}[column sep=8ex]
    C_p \arrow{r}{F} \arrow{d}{d_i} & D_p \arrow{d}{d_i} & D_{p+q+1} \arrow{l}[swap]{d_{p+1}...d_{p+q+1}} \arrow{d}{d_i} \\
    C_{p-1} \arrow{r}{F} & D_{p-1} & D_{p+q} \arrow{l}{d_p...d_{p+q}} 
\end{tikzcd}
\]
It follows that $(F/^hD)_{\bullet,q}$ satisfies the semi-simplicial relation
\begin{center}
    $d_i^{(p)} \circ d_j^{(p)} = d_{j-1}^{(p)} \circ d_i^{(p)}, \,\,\, \text{if} \, \, i < j$.
\end{center}
Similarly in the $q$-direction, for all $0 \leq j \leq q$, the map $d_j^{(q)} : (F/^hD)_{p,q} \rightarrow (F/^hD)_{p,q-1}$ is induced from 
\[
\begin{tikzcd}[column sep=10ex]
    C_p \arrow{r}{F} \arrow{d}{\text{id}} & D_p \arrow{d}{\text{id}} & D_{p+q+1} \arrow{l}[swap]{d_{p+1}...d_{p+q+1}} \arrow{d}{d_{p+1+j}} \\
    C_{p} \arrow{r}{F} & D_{p} & D_{p+q} \arrow{l}{d_{p+1}...d_{p+q}} 
\end{tikzcd}
\]
It can be verified that $(F/^hD)_{p,\bullet}$ satisfies similar relations. This shows that $(F/^hD)_{\bullet,\bullet}$ is a bi-semi-simplicial space. We also denote by $(F/D)_{\bullet,\bullet}$ the bi-semi-simplicial space obtained from taking the strict pullback of diagram \ref{dg 6.11}. 

\begin{lem}   \label{lemma bisemi reso}
    If $F$ is a Reedy fibration or if $D_\bullet$ is Reedy fibrant, then the canonical map 
    \begin{center}
        $(F/D)_{\bullet,\bullet} \rightarrow (F/^hD)_{\bullet,\bullet}$
    \end{center}
    induces a weak homotopy equivalence on their classifying spaces. 
\end{lem}
\proof   
If $F$ is a Reedy fibration or if $D_\bullet$ is Reedy fibrant, then either $F_p$ or $d_{p+1}...d_{p+q+1}$ in diagram \ref{dg 6.11} is a fibration. It follows that the canonical map $(F/D)_{p,q} \rightarrow (F/^hD)_{p,q}$ is a weak homotopy equivalence for all $p, q \geq 0$. Then by \cite[Theorem 2.2]{ebert2019semisimplicial}, the map $\Vert (F/D)_{\bullet,\bullet} \Vert \rightarrow \Vert (F/^hD)_{\bullet,\bullet} \Vert$ is a weak equivalence.

\qed

Note that $(F/^hD)_{\bullet,\bullet}$ has augmentation maps 
\begin{align}
    \epsilon_{p,q} &: (F/^hD)_{p,q} \overset{\text{pr}_1}{\longrightarrow} C_p \\
    \eta_{p,q} &: (F/^hD)_{p,q} \overset{\text{pr}_2}{\longrightarrow} D_{p+q+1} \overset{d_0^{p+1}}{\longrightarrow} D_q
\end{align}
where $\epsilon_{p,q}$ are fibrations for all $p, q \geq 0$. 

\begin{proposition}   \label{prop bi reso weak equi}
    If $D_\bullet$ is a weakly unital semi-Segal space, then the map
    \begin{center}
        $\Vert \epsilon_{\bullet,\bullet} \Vert : \Vert (F/^hD)_{\bullet,\bullet} \Vert \longrightarrow \Vert C_\bullet \Vert$
    \end{center}
    is a weak equivalence. 
\end{proposition}

\proof
We first assume that $D_\bullet$ is Reedy fibrant. Then, together with the Segal condition, it follows that the right square of the diagram
\[
\begin{tikzcd}[column sep=8ex]
    (F/D)_{p,q} \arrow{d}{\widetilde{\epsilon}_{p,q}} \arrow{r} & D_{p+q+1} \arrow{r}{d_0^p} \arrow{d}{d_{p+1}...d_{p+q+1}} & D_{q+1} \arrow{d}{d_1...d_{q+1}} \\
    C_p \arrow{r}{F_p} & D_p \arrow{r}{d_0^p} & D_0 
\end{tikzcd}
\]
is homotopy cartesian. We write $g_q := d_1...d_{q+1}$ and denote by $\widetilde{\epsilon}_{p,q}$ the projection from $(F/D)_{p,q}$ to its first factor. By construction, the left square of the above diagram is also homotopy cartesian and all vertical maps are fibrations. This implies that for $c \in C_p$, the canonical map 
\begin{center}
    $\widetilde{\epsilon}_{p,q}^{-1} (c) \overset{\simeq}{\longrightarrow} g_q^{-1}(F(c_p))$
\end{center}
is a weak homotopy equivalence for all $q \geq 0 $, where we write $c_p:= d_0^p(c) \in C_0$. Since $D_\bullet$ is weakly unital and Reedy fibrant, then in particular there exists a cocartesian 1-simplex $y \in D_1$ with $d_1y = F(c_p)$. We define $y/D_{\bullet+1+1}$ to be the pullback
\[
\begin{tikzcd}
    y/D_{\bullet+1+1} \arrow{r} \arrow{d} & D_{\bullet+1+1} \arrow{d}{d_2...d_{\bullet+2}} \\
    \{y\} \arrow{r} & D_1
\end{tikzcd}
\] 
Note that there is a canonical projection $y/D_{\bullet+1+1} \rightarrow g_{\bullet}^{-1}(F(c_p))$, induced by $d_1$. Then by similar argument as in the proof of \cite[Lemma 5.3]{steimle2021additivity}, the induced map after passing to geometric realizations $\Vert y/D_{\bullet+1+1} \Vert \rightarrow \Vert g_\bullet^{-1}(F(c_p)) \Vert$ is both nullhomotopic and a weak homotopy equivalence. It follows that $\Vert g_\bullet^{-1}(F(c_p)) \Vert$ is weakly contractible, hence so is $ \Vert \widetilde{\epsilon}^{-1}_{p,\bullet} (c) \Vert$, for every $c \in C_p$. Thus by \cite[Lemma 2.14]{ebert2019semisimplicial}, the map $\Vert (F/D)_{p,\bullet} \Vert \rightarrow C_p $ is a weak equivalence at every $p$. Then it follows from \cite[Theorem 2.2]{ebert2019semisimplicial} that the map
\begin{center}
    $\Vert \widetilde{\epsilon}_{\bullet,\bullet} \Vert : \Vert (F/D)_{\bullet,\bullet} \Vert \longrightarrow \Vert C_\bullet \Vert$
\end{center}
is a weak homotopy equivalence. Observe that $\vert \widetilde{\epsilon}_{\bullet,\bullet} \Vert$ factors through $\Vert (F/^hD)_{\bullet,\bullet} \Vert$ via the composition $\Vert (F/D)_{\bullet,\bullet} \Vert \rightarrow \Vert (F/^hD)_{\bullet,\bullet} \Vert  \overset{\Vert \epsilon_{\bullet,\bullet} \Vert}{\longrightarrow} \Vert C_\bullet \Vert$, so by lemma \ref{lemma bisemi reso},  $\Vert \epsilon_{\bullet,\bullet} \Vert$ is also a weak equivalence. 

In the general case, consider a level equivalence $D_\bullet \rightarrow D_\bullet^{\text{Rf}}$ in which $D_\bullet^{\text{Rf}}$ is Reedy fibrant. We denote by $F'$ the composition $C_\bullet \rightarrow D_\bullet \rightarrow D_\bullet^{\text{Rf}}$. By Lemma \ref{lemma stable under reedy}, $D_\bullet^{\text{Rf}}$ is weakly unital, so we can apply the above argument to $D_\bullet^{\text{Rf}}$. Therefore the map $\Vert (F'/^h D^{\text{Rf}})_{\bullet,\bullet} \Vert \rightarrow \Vert C_\bullet \Vert$ is a weak equivalence. Finally we observe that $\Vert \epsilon_{\bullet,\bullet} \Vert$ factors through $\Vert (F'/^h D^{\text{Rf}})_{\bullet,\bullet} \Vert$ via the commutative triangle
\[
\begin{tikzcd}
    & \Vert (F'/^hD^{\text{Rf}})_{\bullet,\bullet} \Vert \arrow{dr}{\simeq} & \\
    \Vert (F/^hD)_{\bullet,\bullet} \Vert \arrow{ur}{\simeq} \arrow{rr}{\Vert \epsilon_{\bullet,\bullet} \Vert } & & \Vert C_\bullet \Vert 
\end{tikzcd}
\]
This proves that $\Vert \epsilon_{\bullet,\bullet} \Vert$ is a weak equivalence. 
\qed

\begin{proposition}   \label{prop homotopy commutative}
    Let $F : C_\bullet \rightarrow D_\bullet$ be a map of semi-simplicial spaces. Then the diagram
    \begin{equation}  \label{dg hocom}
    \begin{tikzcd}
        & \Vert (F/^hD)_{\bullet,\bullet} \Vert \arrow{dl}[swap]{\Vert \epsilon_{\bullet,\bullet} \Vert} \arrow{dr}{\Vert \eta_{\bullet,\bullet} \Vert} & \\
        \Vert C_\bullet \Vert \arrow{rr}{\Vert F\Vert} & & \Vert D_\bullet \Vert 
    \end{tikzcd}
    \end{equation}
    is homotopy commutative. 
\end{proposition}

\proof
Recall that an element of $(F/^hD)_{p,q}$ is a triple $(x, \gamma, y)$, where $x \in C_p, y \in D_{p+q+1}$ and $\gamma$ is a path in $D_p$ such that $\gamma (0) = F(x)$ and $\gamma (1) = d_{p+1}...d_{p+q+1}(y)$. For each $p, q \geq 0$, we define a map
\begin{center}
    $H_{p,q} : I \times (F/^hD)_{p,q} \times \Delta^p \times \Delta^q \longrightarrow \Vert D_\bullet \Vert$
\end{center}
by setting 
\[  H_{p,q}(t;x, \gamma, y ; r ; s) := \left\{
\begin{array}{ll}
      (y ; 2tr , (1-2t)s) \in D_{p+q+1} \times \Delta^{p+q+1} & 0 \leq t \leq \frac{1}{2}  \\
      (\gamma(2-2t) ; r) \in D_p \times \Delta^p & \frac{1}{2} \leq 2 \leq 1
\end{array} \right.
\]
One can check that these maps respect the semi-simplicial relations and hence descend to a homotopy $H : I \times \Vert (F/^hD)_{\bullet,\bullet} \Vert \rightarrow \Vert D_\bullet \Vert$. Moreover, we see that 
\begin{align}
    H_0 &= \Vert \eta_{\bullet,\bullet} \Vert : \Vert (F/^hD)_{\bullet,\bullet} \Vert \rightarrow \Vert D_\bullet \Vert \\
     H_1 &= \Vert F \Vert  \circ \Vert \epsilon_{\bullet,\bullet} \Vert : \Vert (F/^hD)_{\bullet,\bullet} \Vert \rightarrow \Vert C_\bullet \Vert \rightarrow \Vert D_\bullet \Vert 
\end{align}
This shows that $H$ is a homotopy between $\Vert \eta_{\bullet,\bullet} \Vert$ and $\Vert F \Vert \circ \Vert \epsilon_{\bullet,\bullet} \Vert$. Thus diagram \ref{dg hocom} is homotopy commutative. 
\qed

There is one more construction we need before we prove Theorem \ref{thm A} and Theorem \ref{thm B}. Let $D_\bullet$ be a semi-simplicial space. Fix $k \geq 0$. For an object $y \in D_k$, we denote by $D_{\bullet+k+1}/^hy$ the semi-simplicial space in which at each level $D_{q+k+1}/^hy$ is the homotopy fiber product of the diagram
\begin{center}
    $\{y \} \longrightarrow D_k \overset{d_0^{q+1}}{\longleftarrow} D_{q+k+1} $
\end{center}
Similarly, given a semi-simplicial map $F : C_\bullet \rightarrow D_\bullet$, for $y \in D_k$ we define $(F/^hy)_\bullet$ to be the semi-simplicial space whose $p$-th space $(F/^hy)_p$ is the homotopy fiber product of 
\begin{center}
    $C_p \overset{F}{\longrightarrow} D_p \overset{d_{p+1}...d_{p+k+1}}{\xleftarrow{\hspace{1cm}}} D_{p+k+1}/^hy$
\end{center}
We reserve the notation $(F/y)_\bullet$ and $D_{\bullet+k+1}/y$ for the semi-simplicial spaces obtained from taking the strick pullbacks of the same diagrams as above.

\section{Proof of Theorem A}    \label{sec 4}

\textit{Proof of Theorem A.} By Proposition \ref{prop bi reso weak equi} and \ref{prop homotopy commutative}, it is enough to show that the map $\Vert \eta_{\bullet,\bullet} \Vert$ is a weak equivalence. We first assume that $D_\bullet$ is Reedy fibrant and $F$ is a Reedy fibration.  Recall that the augmentation map $\eta_{p,q}$ is the composition
\begin{center}
    $(F/^hD)_{p,q} \overset{\text{pr}_2}{\longrightarrow} D_{p+q+1} \overset{d_0^{p+1}}{\longrightarrow} D_q $ 
\end{center}
so if $D_\bullet$ is Reedy fibrant, then $\eta_{p,q}$ is a fibration. 

Fix $q>0$. Note that by the Segal condition, there exists a homotopy cartesian diagram
\begin{equation}   \label{dg 6.6}
\begin{tikzcd}[column sep =10ex]
    D_{p+q+1} \arrow{d}{d_0^{p+1}} \arrow{r}{d_{p+2}...d_{p+q+1}} & D_{p+1} \arrow{d}{d_0^{p+1}} \\
    D_q \arrow{r}{d_1...d_q} & D_0 
\end{tikzcd}
\end{equation}
Let $y \in D_q$ with $d_1...d_q(y) = y_0  \in D_0$. By construction $D_{p+q+1}/y$ is the actual fiber of the left vertical map  of \ref{dg 6.6} over $y$. Then it follows that the induced map on the vertical fibers
\begin{center}
    $D_{p+q+1}/y \overset{\simeq}{\longrightarrow} D_{p+0+1}/y_0$
\end{center}
is a weak equivalence for all $p \geq 0$. So by \cite[Theorem 2.2]{ebert2019semisimplicial}, the map $\Vert D_{\bullet+q+1}/y \Vert \rightarrow \Vert D_{\bullet+0+1}/y_0 \Vert $ is a weak equivalence. 

On the other hand, by construction $(F/y)_p$ fits into the following diagram
\[
\begin{tikzcd}[column sep=10ex]
    (F/y)_p \arrow{d} \arrow{r} & (F/D)_{p,q} \arrow{d} \arrow{r} & C_p \arrow{d}{F} \\
    D_{p+q+1}/y \arrow{r} \arrow{d} & D_{p+q+1} \arrow{d}{d_0^{p+1}} \arrow{r}{d_{p+1}...d_{p+q+1}} & D_p \\
    \{y\} \arrow{r} & D_q
\end{tikzcd}
\]
in which all smaller squares are strict pullback diagrams, hence so is the left outer square. It follows that $(F/y)_p$ is the fiber of the map $\eta_{p,q} : (F/D)_{p,q} \rightarrow D_q$ over $y$. Then by \cite[Lemma 2.14]{ebert2019semisimplicial}, there exists a homotopy fibration sequence 
\begin{center}
    $\Vert (F/y)_\bullet \Vert \rightarrow \Vert (F/D)_{\bullet,q} \Vert \rightarrow D_q$ 
\end{center}
The proof will be complete if we can show that $\Vert (F/y)_\bullet \Vert$ is weakly contractible, for this would then imply that $q \longmapsto (\Vert (F/D)_{\bullet,q} \Vert  \rightarrow D_q)$ is a level equivalence. To see that $\Vert (F/y)_\bullet \Vert$ is weakly contractible, consider the following diagram
\[
\begin{tikzcd}[column sep=10ex]
    C_p \arrow{d}{\text{id}} \arrow{r}{F} & D_p \arrow{d}{\text{id}} & D_{p+q+1}/y \arrow{d}{\simeq} \arrow{l}[swap]{d_{p+1}...d_{p+q+1}} \\
    C_p \arrow{r}{F} & D_p & D_{p+0+1}/y_0 \arrow{l}[swap]{d_{p+1}} 
\end{tikzcd}
\]
where the right vertical map was shown to be a weak equivalence. This implies that the induced map on pullbacks $(F/y)_p \rightarrow (F/y_0)_p$ is a weak equivalence for every $p$, hence so is the map $\Vert (F/y)_\bullet \Vert \rightarrow \Vert (F/y_0)_\bullet \Vert$. By the assumption that $\Vert (F/^hy_0)_\bullet \Vert$ is contractible, it follows that $\Vert (F/y)_\bullet \Vert$ is weakly contractible. 

In the general case, we can first apply Reedy fibrant replacement to $D_\bullet$, followed by a Reedy fibrant replacement to the resulting composition of maps. This gives rise to a diagram
\[
\begin{tikzcd}[column sep=8ex]
    C_\bullet \arrow{r}{F} \arrow{d}{\simeq} & D_\bullet \arrow{d}{\simeq} \\
    C_\bullet^{\text{Rf}} \arrow{r}{\widetilde{F}} & D_\bullet^{\text{Rf}} 
\end{tikzcd}
\]
in which all vertical maps are level equivalences, $\widetilde{F}$ is a Reedy fibration and $C_{\bullet}^{\text{Rf}}$ and $D_\bullet^{\text{Rf}}$ are Reedy fibrant. It is also easy to see that applying Reedy fibrant replacement does not change the second assumption of Theorem \ref{thm A}, for there is an induced level equivalence $(F/^hy_0)_\bullet \rightarrow (\widetilde{F}/^h\widetilde{y_0})_\bullet$, where $\widetilde{y_0}$ is the image of $y_0$ under the level equivalence. The previous argument applies to $D_\bullet^{\text{Rf}}$, from which we conclude that $\widetilde{F}$ is a weak equivalence. It then follows that $F$ is a weak equivalence.

\qed

\section{Proof of Theorem B}   \label{sec 5}

\textit{Proof of Theorem B}. Since $D_\bullet$ is weakly unital, then by Proposition \ref{prop bi reso weak equi} and \ref{prop homotopy commutative}, it is enough to show that there exists a homotopy fibration sequence 
\begin{center}
    $\Vert (F/^hy_0)_\bullet \Vert \longrightarrow \Vert (F/^hD)_{\bullet,\bullet} \Vert \overset{\Vert \eta_{\bullet,\bullet} \Vert}{\longrightarrow} \Vert D_\bullet \Vert$
\end{center}
for each $y_0 \in D_0$. 

Let $y \in D_q$. Consider the following homotopy commutative cube
\[
\begin{tikzcd}
    & (F/^hy)_p \arrow{dl} \arrow{rr} \arrow{dd} & & D_{p+q+1} /^hy \arrow{dl} \arrow{dd}{d_{p+1}...d_{p+q+1}} \\
    (F/^hD)_{p,q} \arrow{dd} \arrow{rr} & & D_{p+q+1} \arrow{dd}[near start]{d_{p+1}...d_{p+q+1}} & \\
    & C_p \arrow{rr}[near start]{F} \arrow{dl}{\text{id}} & & D_p \arrow{dl}{\text{id}} \\
    C_p \arrow{rr}{F} & & D_p & 
\end{tikzcd}
\]
where the front, back and bottom faces are all homotopy cartesian. Then by the 2-out-of-3 property, so is the top face. On the other hand, we denote by $P_y D_q$ the space of paths in $D_q$ which starts at $y \in D_q$. Observe that the top face also fits into the following commutative diagram
\[
\begin{tikzcd}
    (F/^hy)_p \arrow{r} \arrow{d} & D_{p+q+1}/^hy \arrow{r} \arrow{d} & P_y D_q \arrow{d}{\text{ev}_1} \\
    (F/^hD)_{p,q} \arrow{r} & D_{p+q+1} \arrow{r}{d_0^{p+1}} & D_q
\end{tikzcd}
\]
where $\text{ev}_1 : P_y D_q \rightarrow D_q$ is the map given by $\text{ev}_1(\gamma) = \gamma (1)$. Note that by construction the right square is homotopy cartesian, therefore so is the outer square. Since this is true for all $p \geq 0$, then by \cite[Lemma 2.13]{ebert2019semisimplicial}, together with the fact that $P_y D_q$ is contractible, there exists a homotopy fibration sequence
\begin{equation}   \label{eq fib seq}
    \Vert (F/^hy)_\bullet \Vert \rightarrow \Vert (F/^hD)_{\bullet,q} \Vert \rightarrow D_q
\end{equation}
over $y \in D_q$. 

We next show that the semi-simplicial map $q \mapsto \Vert \eta_{\bullet,q} \Vert$ is homotopy cartesian. We use the criteria in \cite[Lemma 2.11]{ebert2019semisimplicial}. In the first case when $q=i=1$, since diagram \ref{dg 6.6} is homotopy cartesian, so for any $y' \in D_1$, the map $D_{p+1+1}/^hy' \rightarrow D_{p+0+1}/^hd_1y'$ is a weak equivalence for each $p$. This further implies that the induced map
\begin{center}
    $\Vert (F/^hy')_\bullet \Vert \rightarrow \Vert (F/^hd_1y')_\bullet \Vert$
\end{center}
is a weak equivalence. This map fits in to the following commutative cube 
\begin{equation}  \label{dg quillenb} 
\begin{tikzcd}
   & \Vert (F/^hy')_\bullet \Vert \arrow{dd} \arrow{rr} \arrow{dl}{\simeq} & & P_{y'} D_1 \arrow{dl}{d_1} \arrow{dd}{\text{ev}_1} \\
   \Vert (F/^hd_1y')_\bullet \Vert \arrow{dd} \arrow{rr} & & P_{d_1y'} D_0 \arrow[near start]{dd}{\text{ev}_1} & \\
   & \Vert (F/^hD)_{\bullet,1} \Vert \arrow{dl}{d_1} \arrow{rr}[near start]{\Vert \eta_{\bullet,1} \Vert} & & D_1 \arrow{dl}{d_1}  \\
   \Vert (F/^hD)_{\bullet,0} \Vert \arrow{rr}{\Vert \eta_{\bullet,0} \Vert} & & D_0 & 
\end{tikzcd}
\end{equation}
where the front and the back faces are homotopy cartesian. Since $P_{y'} D_1$ and $P_{d_1y'} D_0$ are contractible, then the induced map $\text{hofib}_{y'} (\Vert \eta_{\bullet,1} \Vert) \rightarrow \text{hofib}_{d_1y'} (\Vert \eta_{\bullet,0} \Vert)$ is a weak equivalence. Since this is true for every $y' \in D_1$, thus the bottom square is homotopy cartesian. 

In the second case when $i=0$, we need to show that for any $q > 0$, the square
\[
\begin{tikzcd}[column sep=10ex]
    \Vert (F/^hD)_{\bullet,q} \Vert \arrow{d}{d_0} \arrow{r}{\Vert \eta_{\bullet,q} \Vert} & D_q \arrow{d}{d_0} \\
    \Vert (F/^hD)_{\bullet,q-1} \Vert  \arrow{r}{\Vert \eta_{\bullet,q-1} \Vert} & D_{q-1}
\end{tikzcd}
\]
is homotopy cartesian. By the same reasoning as we did in diagram \ref{dg quillenb}, it is enough to show that for any $y \in D_q$, the map
\begin{center}
    $\Vert (F/^hy)_\bullet \Vert \rightarrow \Vert (F/^hd_0y)_\bullet \Vert$
\end{center}
is a weak equivalence. 

Observe that there exists a commutative diagram 
\begin{equation}   \label{dg last}
\begin{tikzcd}[column sep=10ex]
    \Vert (F/^hy)_\bullet \Vert \arrow{d} \arrow{r} & \Vert (F/^hd_0 y)_\bullet \Vert \arrow{d} \\
    \Vert (F/^hy')_\bullet \Vert \arrow{r}{\simeq} & \Vert (F/^hd_0y')_\bullet \Vert 
\end{tikzcd}
\end{equation}
where we write $y' := d_2...d_q(y)$. The bottom map is a weak equivalence by the second assumption of Theorem \ref{thm B}. We show that all other maps in this diagram are weak equivalences. The left vertical map is induced by the square
\begin{equation}   \label{dg quillenb2}
\begin{tikzcd}[column sep=10ex]
    D_{p+q+1} \arrow{r}{d_{p+3}...d_{p+q+1}} \arrow{d}{d_0^{p+1}} & D_{p+1+1} \arrow{d}{d_0^{p+1}} \\
    D_q \arrow{r}{d_2...d_q} & D_1
\end{tikzcd}
\end{equation}
which is homotopy cartesian by the Segal condition. The induced map on the vertical homotopy fibers $D_{p+q+1}/^hy \rightarrow D_{p+1+1}/^hy'$ is a weak equivalence for every $p$. It follows that $\Vert (F/^hy)_\bullet \Vert \rightarrow \Vert (F/^hy')_\bullet \Vert$ is a weak equivalence. That the right vertical map of diagram \ref{dg last} is a weak equivalence can be proven in the same way. Thus the top map in diagram \ref{dg quillenb2} is a weak equivalence. This completes the proof that the semi-simplicial map $q \mapsto \Vert \eta_{\bullet,q} \Vert$ is homotopy cartesian. 

Therefore, by \cite[Theorem 2.12]{ebert2019semisimplicial}, there exists a homotopy cartesian diagram 
\[
\begin{tikzcd}[column sep=8ex]
    \Vert (F/^hD)_{\bullet,0} \Vert \arrow{r} \arrow{d}{\Vert \eta_{\bullet,0} \Vert} & \Vert (F/^hD)_{\bullet,\bullet} \Vert  \arrow{d}{\Vert \eta_{\bullet,\bullet} \Vert} \\
    D_0 \arrow{r} & \Vert D_\bullet \Vert   
\end{tikzcd}
\]
Finally, from the homotopy fibration sequence in \ref{eq fib seq}, we obtain a homotopy fibration sequence
\begin{center}
    $\Vert (F/^hy_0)_\bullet \Vert \longrightarrow \Vert (F/^hD)_{\bullet,\bullet} \Vert \overset{\Vert \eta_{\bullet,\bullet} \Vert}{\longrightarrow} \Vert D_\bullet \Vert$
\end{center}
for each $y_0 \in D_0$. 

\qed 

As an application, we provide a proof of Corollary \ref{cor c}, which addresses a specific case where the assumptions of Theorem \ref{thm B} are satisfied, thereby enabling the identification of the homotopy fiber. 

\noindent \textit{Proof of Corollary 1.1}. Consider a map of semi-Segal spaces $F: C_\bullet \rightarrow D_\bullet$, where $D_\bullet$ is a weakly unital semi-Segal space in which every 1-simplex is cartesian. We want to show that for any 1-simplex $y \in D_1$,  the induced map 
\begin{equation}  \label{eq map weak} 
    \Vert (F/^hy)_\bullet \Vert \rightarrow \Vert (F/^hd_0y)_\bullet \Vert
\end{equation}
is a weak equivalence. Observe that this map is induced by the following diagram
\[
\begin{tikzcd}[column sep=10ex]
    C_p \arrow{d}{\text{id}} \arrow{r}{F} & D_p \arrow{d}{\text{id}} & D_{p+1+1}/^hy \arrow{l}[swap]{d_{p+1}d_{p+2}} \arrow{d}{d_{p+1}} \\
    C_p \arrow{r}{F} & D_p & D_{p+0+1}/^hd_0y \arrow{l}{d_{p+1}} 
\end{tikzcd}
\]
The proof will be complete if we can show that the right vertical map $D_{p+1+1}/^hy \overset{d_{p+1}}{\longrightarrow} D_{p+0+1}/^hd_0y$ is a weak equivalence for any $y \in D_1$, for this would imply that $(F/^hy)_\bullet \rightarrow (F/^hd_0y)_\bullet$ is a level equivalence. But this is equivalent to proving that the square
\[
\begin{tikzcd}[column sep=10ex]
    D_{p+2} \arrow{d}{d_0^{p+1}} \arrow{r}{d_{p+1}} & D_{p+1} \arrow{d}{d_0^{p+1}} \\
    D_1 \arrow{r}{d_0} & D_0
\end{tikzcd}
\]
is homotopy cartesian. Observe that the above square fits into the following larger diagram
\[
\begin{tikzcd}[column sep=1ex, row sep=1ex]
    D_{p+2} \arrow{rr}{d_{p+1}} \arrow{dd}{d_0^p} & & D_{p+1} \arrow{dd}{d_0^p} \arrow{rr}{d_{p+1}} &  & D_p \arrow{dd}{d_0^p} \\
    & (A) & & (B) & \\
    D_2 \arrow{dd}{d_0} \arrow{rr}{d_1} &  & D_1 \arrow{dd}{d_0} \arrow{rr}{d_1} &  & D_0 \\
    & (C) & & & \\
    D_1 \arrow{rr}{d_0} &  & D_0 & & 
\end{tikzcd}
\]
Here the squares $(B)$ and $(A+B)$ are homotopy cartesian by the Segal condition, so by \cite[Lemma 3.3]{steinebrunner2022locally}, so is the square $(A)$. By the assumption that every 1-simplex in $D_1$ is cartesian, the square $(C)$ is homotopy cartesian, hence so is $(A+C)$. 

\qed

\bibliography{ms}
\bibliographystyle{amsalpha}

\bigskip
{\footnotesize
Yuxun Sun, \textsc{Department of Mathematics, University at Buffalo, Buffalo, NY, 14260} \par\nopagebreak
  \textit{E-mail address}:  \texttt{yuxunsun@buffalo.edu}.
}

\end{document}